\documentclass{amsart}

\usepackage{fancyhdr}
\usepackage{lastpage}
\usepackage{stmaryrd,yhmath}

\newcommand{\bdis}{\begin{displaymath}}
\newcommand{\edis}{\end{displaymath}}
\newcommand{\be}{\begin{equation}}
\newcommand{\ee}{\end{equation}}

\newcommand{\mcal}{\mathcal}

\newcommand{\vp}{\varphi}

\newcommand{\mT}{\mathring{T}}

\newcommand{\zf}{\zeta\left(\frac{1}{2}+it\right)}
\newcommand{\zfvp}{\zeta\left(\frac{1}{2}+i\vp_1(t)\right)}

\newtheorem{theorem}{Theorem}
\newtheorem{lemma}[]{Lemma}

\theoremstyle{definition}

\newtheorem{cor}[]{Corollary}

\theoremstyle{remark}
\newtheorem{remark}[]{Remark}

\newtheorem*{mydef2}{{\bf Definition}}

\numberwithin{equation}{section}

\begin{document}

\title{Jacob's ladders, the structure of the Hardy-Littlewood integral and some new class of nonlinear integral equations}

\author{Jan Moser}

\address{Department of Mathematical Analysis and Numerical Mathematics, Comenius University, Mlynska Dolina M105, 842 48 Bratislava, SLOVAKIA}

\email{jan.mozer@fmph.uniba.sk}


\begin{abstract}
In this paper we obtain new formulae for short and microscopic parts of the Hardy-Littlewood integral, and the first asymptotic formula for the sixth order
expression $|\zeta(\frac{1}{2}+i\vp_1(t))|^4|\zf|^2$. These formulae cannot be obtained in the theories of Balasubramanian, Heath-Brown and Ivic. \\

Dedicated to the 75th aniversary of Anatolii Alekseevich Karatsuba.
\end{abstract}

\maketitle

\section{Introduction}

Let us remind that Hardy and Littlewood started to study the following integral in 1918

\be \label{1.1}
\int_0^T\left|\zf\right|^2{\rm d}t=\int_0^TZ^2(t){\rm d}t
\ee

where

\be\label{1.2}
Z(t)=e^{i\vartheta(t)}\zf,\ \vartheta(t)=-\frac t2\ln\pi+\text{Im}\ln\Gamma\left(\frac 14+i\frac t2\right) ,
\ee

and they have derived the following formula (\cite{1}, pp. 122, 151-156)

\be \label{1.3}
\int_0^T Z^2(t){\rm d}t\sim T\ln T,\ T\to\infty .
\ee

We have shown in the paper \cite{7} that except the asymptotic formula (\ref{1.3}) possessing an unbounded error there is an infinite family of other
asymptotic representations of the Hardy-Littlewood integral (\ref{1.1}). Each member of this family is an almost exact representation of the integral
(\ref{1.1}). Namely, the following formula

\bdis
\int_0^T Z^2(t){\rm d}t=\frac{\vp(T)}{2}\ln\frac{\vp(T)}{2}+(c-\ln 2\pi)\frac{\vp(T)}{2}+c_0+\mcal{O}\left(\frac{\ln T}{T}\right)
\edis

takes place where $\vp(T)$ is the Jacob's ladder, i.e. an arbitrary solution to the nonlinear integral equation

\bdis
\int_0^{\mu[x(T)]}Z^2(t)e^{-\frac{2}{x(T)}t}{\rm d}t=\int_0^T Z^2(t){\rm d}t ,
\edis

and $\mu(y)\geq 7y\ln y;\ \mu(y)\to y=\vp_\mu(T)=\vp(T)$. \\

We obtain new properties of the signal (\ref{1.2}) generated by the Riemann zeta-function. Namely:

\begin{itemize}
\item[(A)] We obtain the multiplicative asymptotic formula

\bdis
\int_T^{T+U}Z^2(t){\rm d}t\sim U\ln T\tan[\alpha(T,U)],\ U\in \left(\left. 0,\frac{T}{\ln T}\right]\right.,\ T\to\infty ,
\edis

in the parts 2-7   of this paper. The application on microscopic ($0<U<1$) and short ($1\leq U<T^{1/3+2\epsilon}$) parts of the Hardy-Littlewood integral
(\ref{1.1}) is the main aim of this formula.

\item[(B)] We also obtain, in the parts 8 - 10 of this work, the formula

\bdis
\int_T^{T+U_1}\left|\zfvp\right|^4\left|\zf\right|^2{\rm d}t\sim \frac{1}{2\pi}U_1\ln^5T,\ T\to\infty
\edis

where $U_1=T^{7/8+2\epsilon},\ \vp_1(t)=\frac 12\vp(t)$. This formula is the first integral asymptotic formula in the theory of the Riemann zeta-function
for the sixth-order expression $\left|\zfvp\right|^4\left|\zf\right|^2$.

\item[(C)] In the part 11 of this work the following property, for example,
is noticed: the Jacob's ladder $\vp_1$ is the asymptotic solution of the nonlinear integral equation

\bdis
\int_{x^{-1}(T)}^{x^{-1}(T+2)}\frac{[\mathcal{T}_n(x(t)-T-1)]^2}{\sqrt{1-(x(t)-T-1)^2}}\left|\zf\right|^2{\rm d}t=\frac{\pi}{2}\ln T ,
\edis

$n=1,2,\dots $, where $\mathcal{T}_n(t),\ t\in [-1,1]$ is the Chebyshev polynomial of the first kind, i.e. the following asymptotic formula

\bdis
\int_{\vp_1^{-1}(T)}^{\vp_1^{-1}(T+2)}\frac{[\mcal{T}_n(\vp_1(t)-T-1)]^2}{\sqrt{1-(\vp_1(t)-T-1)^2}}\left|\zf\right|^2{\rm d}t\sim\frac{\pi}{2}\ln T,\ T\to\infty
\edis
holds true.

\end{itemize}

\section{Necessity of a new expression for the short and the microscopic parts of the Hardy-Littlewood integral}

The Balasubramanian's formula

\bdis
\int_0^T Z^2(t){\rm d}t=T\ln T+(2c-1-\ln 2\pi)T+\mcal{O}(T^{1/3+\epsilon})
\edis

implies (comp. \cite{7}, (2.5), (8.3))

\be \label{2.1}
\int_T^{T+U_0} Z^2(t){\rm d}t=U_0\ln T+ (2c-\ln 2\pi)U_0+\mcal{O}(T^{1/3+\epsilon}),\ U_0=T^{1/3+2\epsilon}
\ee

where $c$ is the Euler's constant. Furthermore, let us remind the Heath-Brown's estimate (see \cite{4}, (7.20), p. 178)

\be\label{2.2}
\begin{split}
& \int_{T-G}^{T+G} Z^2(t){\rm d}t=  \\
& \mcal{O}\left\{ G\ln T+G\sum_K (TK)^{-1/2}\left(S(K)+K^{-1}\int_0^K|S(x)|{\rm d}x\right)e^{-G^2K/T}\right\}
\end{split}
\ee
(for the definition of used symbols see \cite{4}, (7.21)-(7.23)), uniformly for $T^\epsilon\leq G\leq T^{1/2+\epsilon}$.
And, finally, we add the Ivic' estimate (\cite{4}, (7.62))

\be \label{2.3}
\int_{T-G}^{T+G} Z^2(t){\rm d}t=\mcal{O}(G\ln^2T),\ G\geq T^{1/3-\epsilon_0},\ \epsilon_0=\frac{1}{108}\approx 0.009 .
\ee

\begin{remark}
It is obvious that the short intervals $[T-G,T+G]$ with $G=1000$, for example, are not included in the methods od Balasubramanian, Heath-Brown and
Ivic, as well, leading to (\ref{2.1})-(\ref{2.3}).
\end{remark}

In this paper we present a new method how to deal with short and microscopic parts

\be \label{2.4}
\int_T^{T+U} Z^2(t){\rm d}t
\ee

of the Hardy-Littlewood integral (\ref{1.1}). In order to attain this goal we will use only elementary geometric properties of the Jacob's ladders. The basic
idea is expressed in the following theorem.

\begin{theorem}
For $\mu[\vp]=a\vp\ln\vp,\ a\in [7,8]$ the following is true
\be \label{2.5}
\begin{split}
& \int_T^{T+U} Z^2(t){\rm d}t=\left\{ 1+\mcal{O}\left(\frac{\ln\ln T}{\ln T}\right)\right\}U\ln T\tan[\alpha(T,U)], \\
& U\in \left(\left. 0,\frac{T}{\ln T}\right.\right],\ T\geq T_0[\vp]
\end{split}
\ee
where $\alpha=\alpha(T,U)$ is the angle of the chord of the curve $y=\frac 12\vp(T)$ that binds the points $[T,\frac 12\vp(T)]$ and
$[T+U,\frac 12\vp(T+U)]$.
\end{theorem}

\section{A geometric criterion for validity of the usual mean-value theorem}

\subsection{}

First of all we will show the canonical equivalence that follows from (\ref{2.5}). Let us remind (see \cite{7}, (8.3)) that we call the chord binding
the points

\be \label{3.1}
\left[T,\frac 12\vp(T)\right];\ \left[T+U_0,\frac 12\vp(T+U_0)\right]
\ee

of the Jacob's ladder $y=\frac 12\vp(T)$ \emph{the fundamental chord}. By comparison of the formulae (\ref{2.1}) and (\ref{2.5}), $U=U_0$, we obtain
the following asymptotic formula

\be \label{3.2}
\tan[\alpha(T,U_0)]=\frac{\vp(T+U_0)-\vp(T)}{2U_0}=1+\mcal{O}\left(\frac{\ln\ln T}{\ln T}\right) .
\ee

\begin{mydef2}
The chord binding the points
\be \label{3.3}
\left[N,\frac 12\vp(N)\right];\ \left[M,\frac 12\vp(M)\right],\ [N,M]\subset [T,T+U_0]
\ee
which fulfills the property
\bdis
\tan[\alpha(N,M-N)]=1+o(1),\ T\to\infty
\edis
(compare (\ref{3.2})) is called \emph{the almost parallel chord to the fundamental chord}. This property will be denoted by the symbol $\fatslash$.
\end{mydef2}

\begin{cor}
Let $[N,M]\subset [T,T+U_0]$. Then
\be \label{3.4}
\frac{1}{M-N}\int_N^M Z^2(t){\rm d}t \sim \ln T \quad \Leftrightarrow \quad \fatslash .
\ee
\end{cor}

\begin{remark}
Wee see that the analytic property
\bdis
\frac{1}{M-N}\int_N^M Z^2(t){\rm d}t \sim \ln T
\edis
(the usual mean-value theorem) is equivalent to the geometric property $\fatslash$ of the Jacob's ladder $y=\frac 12\vp(T)$.
\end{remark}

\subsection{}

Let us consider the set of all chords of the curve $y=\frac 12\vp(T)$ which are almost parallel to the fundamental chord. Let the generic chord of this
set bind the points (\ref{3.3}). Then, from Corollary 1, we obtain

\begin{cor}
There is an continuum of intervals $[N,M]\subset [T,T+U_0]$ such that the following asymptotic formula
\be \label{3.5}
\int_N^M Z^2(t){\rm d}t\sim (M-N)\ln T
\ee
holds true.
\end{cor}

\begin{remark}
There is, for example, continuum of intervals $[N,M]:\ 0<M-N<1$ such that the asymptotic formula (\ref{3.5}) holds true (chords approaching to zero).
\end{remark}

\section{On the microscopic parts of the Hardy-Littlewood integral in the neighbourhoods of zeroes of the function $\zf$}

Let $\gamma,\gamma'$ be a pair of neighbouring zeroes of the function $\zf$. The function $\frac{1}{2}\vp(T)$ is necessarily convex on some right
neighbourhood of the point $T=\gamma$, and this function is necessarily concave on some left neighbourhood of the point $T=\gamma'$. Therefore, there
exists the minimal value $\rho\in (\gamma,\gamma')$ such that $[\rho,\frac 12\vp(\rho)]$ is the inflexion point of the curve $y=\frac 12\vp(t)$. At this
point, by properties of the Jacob's ladders, we have $\vp'(\rho)>0$. Let furthermore $\beta=\beta(\gamma,\rho)$ be the angle of the chord binding
the points

\be \label{4.1}
\left[\gamma,\frac 12\vp(\gamma)\right];\  \left[\rho,\frac 12\vp(\rho)\right] .
\ee

Then we obtain from Theorem 1

\begin{cor}
For every sufficiently big zero $T=\gamma$ of the function $\zeta(1+\frac 12 iT)$ the following formulae describing the microscopic parts of the
Hardy-Littlewood integral hold true
\begin{itemize}

\item[(A)] continuum of formulae
\be \label{4.2}
\int_\gamma^{\gamma+U} Z^2(t){\rm d}t\sim U\ln \gamma \tan\alpha ,\ \alpha \in (0,\beta(\gamma,\rho)),\ U=U(\gamma,\alpha)\in (0,\rho-\gamma)
\ee
where $\alpha$ is the angle of the rotating chord binding the points $[\gamma,\frac 12\vp(\gamma)]$, $[\gamma+U,\frac 12\vp(\gamma+U)]$,

\item[(B)] continuum of the formulae for the chords parallel to the chord given by the points (\ref{4.1})
\be \label{4.3}
\int_N^M Z^2(t){\rm d}t\sim (M-N)\ln N\tan[\beta(\gamma,\rho)],\ \gamma\leq N<M<\gamma'  .
\ee
\end{itemize}
\end{cor}

\begin{remark}
The notion \emph{microscopic parts} of the Hardy-Littlewood integral has its natural origin in the following: by Karatsuba's selbergian estimate (see \cite{5},
p. 265) for \emph{almost all} intervals $[\gamma,\gamma']\subset [T,T+T^{1/3+2\epsilon}]$ we have
\be \label{4.4}
\gamma'-\gamma<A\frac{\ln\ln T}{\ln T}\to 0,\ T\to\infty .
\ee
\end{remark}

\begin{remark}
In connection with (\ref{4.4}) we can remind that if the Riemann hypothesis is true then the Littlewood's estimate takes place (see \cite{6})
\bdis
\gamma'-\gamma < \frac{A}{\ln\ln\gamma},\ \gamma\to\infty .
\edis
\end{remark}

\section{Second class of formulae for parts of the Hardy-Littlewood integral beginning in zeroes of the function $\zf$}

Let $\gamma,\bar{\gamma}$ be a pair of zeroes of the function $\zf$ such that $\bar{\gamma}$ obeys the following conditions

\bdis
\bar{\gamma}=\gamma+\gamma^{1/3+2\epsilon}+\Delta(\gamma),\ 0\leq \Delta(\gamma)=\mcal{O}(\gamma^{1/4+\epsilon}) ,
\edis

(it is sufficient to use the classical Hardy-Littlewood's estimate for the distance between the neighbouring zeroes, \cite{1}, pp. 125, 177-184).
Consequently

\be \label{5.1}
U(\gamma)=\gamma^{1/3+2\epsilon}+\Delta(\gamma)\sim \gamma^{1/3+2\epsilon},\ \gamma\to\infty .
\ee

For the chord binding the points

\be \label{5.2}
\left[\gamma,\frac 12\vp(\gamma)\right],\ \left[\bar{\gamma},\frac 12\vp(\bar{\gamma})\right]
\ee

we have by (\ref{3.1}), (\ref{5.1})

\be \label{5.3}
\tan[\alpha(\gamma,U(\gamma))]=1+\mcal{O}\left(\frac{1}{\ln\gamma}\right) .
\ee

The continuous curve $y=\frac 12\vp(T)$ lies bellow the chord given by the points (\ref{5.2}) on some right neighbourhood of the point $T=\gamma$, and this
curve lies above that chord on some left neighbourhood of the point $T=\bar{\gamma}$. Therefore, there exists a common point of the curve and of the chord.
Let $[\bar{\rho},\frac 12\vp(\bar{\rho})],\ \bar{\rho}\in (\gamma,\bar{\gamma})$ be such a common point that is the closest one to the point
$[\gamma,\frac 12\vp(\gamma)]$. Then we obtain from Theorem 1 the following

\begin{cor}
For every sufficient big zero $T=\gamma$ of the function $\zeta(\frac 12+iT)$ we have the following formulae for the parts (\ref{2.4}) of the
Hardy-Littlewood integral (\ref{1.1})
\begin{itemize}
\item[(A)] continuum of formulae for the rotating chord
\be \label{5.4}
\int_\gamma^{\gamma+U} Z^2(t){\rm d}t\sim U\ln\gamma\tan\alpha ,\  \alpha\in [\eta,1-\eta],\ U=U(\gamma,\alpha)\in (0,\bar{\rho}-\gamma)
\ee
where $\alpha=\alpha(\gamma,U)$ is the angle of the rotating chord binding the points $[\gamma,\frac 12\vp(\gamma)]$ and
$[\gamma+U,\frac 12\vp(\gamma+U)]$, and $0<\eta$ is an arbitrarily small number,

\item[(B)] continuum of formulae for the chords parallel (and almost parallel) to the chord binding the points (\ref{5.2})
\be \label{5.5}
\int_N^M Z^2(t){\rm d}t\sim (M-N)\ln N,\ \gamma\leq N<M\leq \bar{\gamma} .
\ee
\end{itemize}
\end{cor}

\begin{remark}
For example, in the case $\alpha=\pi/6$ we have from (\ref{5.4})
\bdis
\int_\gamma^{\gamma+U(\gamma,\pi/6)} Z^2(t){\rm d}t\sim \frac{1}{\sqrt{3}}U\ln\gamma
\edis
for every sufficiently big zero $T=\gamma$ of the function $\zeta(\frac 12+iT)$.
\end{remark}

\begin{remark}
It is clear that the asymptotic formulae (\ref{2.5}), (\ref{3.4}), (\ref{4.2}), (\ref{4.3}), (\ref{5.4}) and (\ref{5.5}) cannot be derived within
complicated methods of Balasubramanian, Heath-Brown and Ivic.
\end{remark}

\section{An estimate for $\Phi''_{\vp\vp}[\vp(T)]$}

Let us remind (see \cite{7}, (3.5), (3.9)) that

\be\label{6.1}
Z^2(t)=\Phi'_\vp[\vp(t)]\frac{{\rm d}\vp(t)}{{\rm d}t} ,
\ee

where

\be \label{6.2}
\begin{split}
& \Phi'_\vp[\vp]=\frac{2}{\vp^2}\int_0^{\mu[\vp]} te^{-\frac{2}{\vp}t}Z^2(t){\rm d}t
+ Z^2\{\mu[\vp]\}e^{-\frac{2}{\vp}\mu[\vp]}\frac{{\rm d}\mu(\vp)}{{\rm d}\vp} .
\end{split}
\ee

The following lemma is true.

\begin{lemma}
If $\mu[\vp]=a\vp\ln\vp,\ a\in [7,8]$ then
\be \label{6.3}
\Phi''_{\vp\vp}[\vp(T)]=\mcal{O}\left(\frac{1}{\vp}\ln\vp\ln\ln\vp\right),\ T\geq T_0[\vp]
\ee
uniformly with respect to $a$.
\end{lemma}

\begin{remark}
The segment $[7,8]$ is sufficient to our purpose since the continuum of Jacob's ladders corresponds to this segment.
\end{remark}

\begin{proof}
First of all we have (see (\ref{6.2}))
\be \label{6.4}
\Phi''_{\vp\vp}[\vp]=\frac{4}{\vp^3}\int_0^{\mu[\vp]}t\left(\frac{t}{\vp}-1\right)e^{-\frac{2}{\vp}t}Z^2(t){\rm d}t+Q[\vp] ,
\ee

\be\label{6.5}
\begin{split}
& Q[\vp]=e^{-\frac{2}{\vp}\mu[\vp]} \left\{\frac{4}{\vp^2}Z^2\{\mu[\vp]\}\mu[\vp]\frac{{\rm d}\mu[\vp]}{{\rm d}\vp}-\frac{2}{\vp}
Z^2\{\mu[\vp]\}\left(\frac{{\rm d}\mu[\vp]}{{\rm d}\vp}\right)^2+\right. \\
& \left. 2Z\{\mu[\vp]\}Z'_\mu\{\mu[\vp]\}\left(\frac{{\rm d}\mu[\vp]}{{\rm d}\vp}\right)^2+Z^2\{\mu[\vp]\}\frac{{\rm d^2}\mu[\vp]}{{\rm d}\vp^2}\right\} .
\end{split}
\ee
Let
\bdis
g(t)=t\left(\frac{t}{\vp}-1\right)e^{-\frac{2}{\vp}t},\ t\in [0,\mu[\vp]] .
\edis
We apply the following elementary facts
\be \label{6.6}
\begin{split}
& g(0)=g(\vp)=0,\ g'\left[\left( 1-\frac{1}{\sqrt{2}}\right)\vp\right]=g'\left[\left( 1+\frac{1}{\sqrt{2}}\right)\vp\right]=0, \\
& \min\{ g(t)\}=-\frac{1}{\sqrt{2}}\left( 1-\frac{1}{\sqrt{2}}\right)e^{-2+\sqrt{2}}\vp, \\
& \max\{ g(t)\}=\frac{1}{\sqrt{2}}\left( 1+\frac{1}{\sqrt{2}}\right)e^{-2-\sqrt{2}}\vp, \\
& g(t)\leq g(\vp\ln\ln\vp)<\vp\left(\frac{\ln\ln\vp}{\ln\vp}\right)^2,\ t\in [\vp\ln\ln\vp,8\vp\ln\vp] , \\
& Z(t),Z'(t)=\mcal{O}(t^{1/4}) ,
\end{split}
\ee
and the Hardy-Littlewood formula (\ref{1.3}). We have
\be \label{6.7}
\begin{split}
& \frac{4}{\vp^3}\int_0^{\vp\ln\ln\vp}=\mcal{O}\left(\frac{1}{\vp^2}\int_0^{\vp\ln\ln\vp} Z^2(t){\rm d}t\right)=
\mcal{O}\left(\frac{1}{\vp}\ln\vp\ln\ln\vp\right) , \\
& \frac{4}{\vp^3}\int_{\vp\ln\ln\vp}^{8\vp\ln\vp}=\mcal{O}
\left\{\frac{1}{\vp^3}\vp\left(\frac{\ln\ln\vp}{\ln\vp}\right)^2\vp\ln^2\vp\right\}=\mcal{O}\left\{\frac{1}{\vp}(\ln\ln\vp)^2\right\}
\end{split}
\ee
by (\ref{1.3}), (\ref{6.6}) and (see (\ref{6.5}))
\be \label{6.8}
Q[\vp]=\mcal{O}(\vp^{-13})\to 0,\ T\to\infty .
\ee
Finally, we obtain (\ref{6.3}) from (\ref{6.4}) by (\ref{6.7}), (\ref{6.8}).
\end{proof}

\section{Proof of Theorem 1}

By (\ref{6.1}) we have

\bdis
\int_T^{T+U} Z^2(t){\rm d}t=\Phi'_\vp[\vp(t_1)]\int_T^{T+U}{\rm d}\vp=\Phi'_\vp[\vp(t_1)]\{\vp(T+U)-\vp(T)\},
\edis

i.e.

\be \label{7.1}
\begin{split}
& \int_T^{T+U} Z^2(t){\rm d}t=2U\Phi'_\vp[\vp(t_1)]\tan[\alpha(T,U)],\ t_1=t_1(U)\in (T,T+U) , \\
& \tan[\alpha(T,U)]=\frac{\vp(T+U)-\vp(T)}{2U} .
\end{split}
\ee

Next, we have

\be \label{7.2}
\int_T^{T+U_0} Z^2(t){\rm d}t=2U_0\Phi_\vp[\vp(t_2)]\left\{ 1+\mcal{O}\left(\frac{1}{\ln T}\right)\right\},\ t_2=t_2(U_0)\in (T,T+U_0) ,
\ee

by (\ref{3.2}), (\ref{7.1}). Hence, by the comparison of the formulae (\ref{2.1}) and (\ref{7.2}) we obtain

\be \label{7.3}
\Phi'_\vp[\vp(t_2)]=\frac 12\ln T+\mcal{O}(1) .
\ee

Next, from the formula (see \cite{7}, (6.2))

\be \label{7.4}
T-\frac{\vp(T)}{2}\sim (1-c)\pi(T);\ T\sim \frac{\vp(T)}{2} ,
\ee

we obtain

\be \label{7.5}
\vp(t_1)-\vp(t_2)=2(t_1-t_2)+\mcal{O}\left(\frac{T}{\ln T}\right)=\mcal{O}\left(\frac{T}{\ln T}\right), \ U\in\left(\left. 0,\frac{T}{\ln T}\right.\right] ,
\ee

and subsequently (see (\ref{6.3})

\be\label{7.6}
\Phi'_\vp[\vp(t_1)]-\Phi'_\vp[\vp(t_2)]=\mcal{O}\{|\Phi''_{\vp\vp}(T)|\cdot |\vp(t_1)-\vp(t_2)|\}=\mcal{O}(\ln\ln T) .
\ee

Therefore we obtain

\be \label{7.7}
\Phi'_\vp[\vp(t_1)]=\frac 12\ln T+\mcal{O}(\ln\ln T) ,
\ee

by (\ref{7.3}), (\ref{7.6}). Finally, (\ref{2.5}) follows from (\ref{7.1}), (\ref{7.7}).

\begin{remark}
Similarly to (\ref{7.6}) we have
\be \label{7.8}
\Phi'_\vp[\vp(t_1)]-\Phi'_\vp[\vp(t)]=\mcal{O}(\ln\ln T),\ t\in [T,T+U] ,
\ee
and obtain (see (\ref{6.1}), (\ref{7.7}), (\ref{7.8}))
\be \label{7.9}
Z^2(t)=\frac 12\left\{ 1+\mcal{O}\left(\frac{\ln\ln t}{\ln t}\right)\right\}\ln t\frac{{\rm d}\vp(t)}{{\rm d}t},\ t\in [T,T+U],\
U\in\left(\left. 0,\frac{T}{\ln T}\right.\right] .
\ee
\end{remark}

\section{The integral asymptotic formula that contains the expression of the sixth order $|\zeta(\frac 12+i\vp_1(t))|^4|\zf|^2$}

\subsection{}

Let us remind that Hardy and Littlewood started to study the following integral in 1926

\bdis
\int_1^T \left|\zf\right|^4{\rm d}t=\int_1^T Z^4(t){\rm d}t ,
\edis

and they derived the following estimate (see \cite{2}, pp. 41,59; \cite{9}, p. 124)

\bdis
\int_1^T\left|\zf\right|^4{\rm d}t=\mcal{O}(T\ln^4T) .
\edis

In 1926 Ingham derived the asymptotic formula

\be \label{8.1}
\int_1^T\left|\zf\right|^4{\rm d}t=\frac{1}{2\pi^2}T\ln^4T+\mcal{O}(T\ln^3T)
\ee

(see \cite{3}, p. 277, \cite{9}, p. 129). Let us remind, finally, the Ingham - Heath-Brown formula (see \cite{4}, p. 129)

\be \label{8.2}
\int_0^T Z^4(t){\rm d}t=T\sum_{k=0}^4 C_k\ln^{4-k}T+\mcal{O}(T^{7/8+\epsilon}),\ C_0=\frac{1}{2\pi^2} ,
\ee

which improves the Ingham formula (\ref{8.1}), (the small improvements of the exponents $1/3$ and $7/8$, see (\ref{2.1}), (\ref{8.2}) are irrelevant
for our purpose).

\subsection{}

In this direction, the following theorem holds true.

\begin{theorem}
\be \label{8.3}
\begin{split}
& \int_T^{T+U_1}\left|\zeta\left(\frac 12+i\vp_1(t)\right)\right|^4\left|\zf\right|^2{\rm d}t\sim \frac{1}{2\pi^2}U_1\ln^5T, \\
& U_1=T^{7/8+2\epsilon},\ \vp_1(t)=\frac 12\vp(t),\ T\to\infty ,
\end{split}
\ee

and the distance of the interaction of the functions

\bdis
\left|\zeta\left(\frac 12+i\vp_1(t)\right)\right|^4,\ \left|\zf\right|^2
\edis

is

\be \label{8.4}
t-\vp_1(t)\sim (1-c)\pi(t) ,
\ee

$c$ is the Euler's constant and $\pi(t)$ is the prime-counting function.
\end{theorem}

\begin{remark}
The formula (\ref{8.3}) is the first integral asymptotic formula in the theory of the Riemann zeta-function for the sixth-order expression
$\left|\zeta\left(\frac 12+i\vp_1(t)\right)\right|^4\left|\zf\right|^2$. This formula cannot be obtained by the methods of Balasubramanian, Heath-Brown and
Ivic.
\end{remark}

\subsection{}

Since (see (\ref{8.4}))

\bdis
T+U_1-\vp_1(T+U_1)\sim (1-c)\pi(T+U_1),\ U_1=T^{7/8+2\epsilon} ,
\edis

we obtain

\bdis
T-\vp_1(T+U_1)\sim (1-c)\pi(T+U_1)-U_1\sim (1-c)\pi(T) ,
\edis

and consequently

\be \label{8.5}
\rho\{[\vp_1(T),\vp_1(T+U)];[T,T+U]\}\sim (1-c)\pi(T);\ \vp_1(T+U_1)<T
\ee

where $\rho$ denotes the distance of the corresponding segments. Next, by using the mean-value theorem in (\ref{8.3}), we obtain

\begin{cor}
\be \label{8.6}
\left|\zeta\left(\frac 12+i\vp_1(\omega)\right)\right|^4\left|\zeta\left(\frac 12+i\omega\right)\right|^2\sim \frac{1}{2\pi^2}\ln^5T
\ee

where

\bdis
\omega\in (T,T+U_1), \vp_1(\omega)\in (\vp_1(T),\vp_1(T+U_1)), \ \omega=\omega(T,U_1,\vp_1) .
\edis
\end{cor}

\begin{remark}
Some \emph{nonlocal interaction} of the functions
\bdis
\left|\zeta\left(\frac 12+i\vp_1(t)\right)\right|^4,\ \left|\zf\right|^2
\edis
is expressed by the formula (\ref{8.6}). This interaction is connected with two segments unboundedly receding each from other (see (\ref{8.5});
$\rho\to\infty$ as $T\to\infty$) - like mutually receding galaxies in the Friedman's expanding Universe.
\end{remark}

\begin{remark}
Since $T\sim\omega,\ \omega\in (T,T+U_1)$ then from (\ref{8.6}) we obtain
\be \label{8.7}
\left|\zeta\left(\frac 12+i\omega\right)\right|\sim\frac{1}{\sqrt{2}\pi}\frac{\ln^{5/2}\omega}{\left|\zeta\left(\frac 12+i\vp_1(\omega)\right)\right|} ,
\ee
i.e. we have \emph{the prediction} of the values $|\zeta(1/2+i\omega)|,\ \omega\in(T,T+U)$ by means of the values $|\zeta(1/2+i\vp_1(\omega))|$
corresponding to the argument $\vp_1(\omega)\in(\vp_1(T),\vp_1(T+U))$ which descend from the very deep past (see (\ref{8.5}), (\ref{8.7})) and - vice
versa.
\end{remark}

\section{Contact point of $|\zf|^2$ with the class of the L-integrable functions of the constant sign}

Let

\be \label{9.1}
\tilde{Z}^2(t)=\frac{{\rm d}\vp_1(t)}{{\rm d}t},\ \vp_1(t)=\frac 12\vp(t),\ t\geq T_0[\vp] ,
\ee

where

\be \label{9.2}
\begin{split}
& \tilde{Z}^2(t)=\frac{Z^2(t)}{2\Phi'_\vp[\vp(t)]}=\frac{\left|\zf\right|^2}{\left\{ 1+\mcal{O}\left(\frac{\ln\ln t}{\ln t}\right)\right\}\ln t} , \\
& t\in [T,T+U],\ U\in \left(\left. 0,\frac{T}{\ln T}\right.\right]
\end{split}
\ee

(see (\ref{6.1}), (\ref{7.7}), (\ref{7.8})). The following lemma holds true (see (\ref{9.1}))

\begin{lemma}
For every integrable function (in the Lebesgue sense) $f(x),\ x\in [\vp_1(T),\vp_1(T+U)]$ the following is true

\be \label{9.3}
\int_T^{T+U}f[\vp_1(t)]\tilde{Z}^2(t){\rm d}t=\int_{\vp_1(T)}^{\vp_1(T+U)}f(x){\rm d}x,\ T\geq T_0[\vp]
\ee

where

\be \label{9.4}
t-\vp_1(t)\sim (1-c)\pi(t) .
\ee
\end{lemma}

\begin{remark}
The formula (\ref{9.3}) is true also in the case of relatively convergent improper Riemann' integral on its right-hand side.
\end{remark}

If $\vp_1\{[\mT,\widering{T+U}]\}=[T,T+U]$ then we have the following formula (see (\ref{9.3}))

\begin{lemma}
\be \label{9.5}
\int_{\mT}^{\widering{T+U}}f[\vp_1(t)]\tilde{Z}^2(t){\rm d}t=\int_T^{T+U}f(x){\rm d}x,\ T\geq T_0[\vp] .
\ee
\end{lemma}

Next, the following lemma holds true.

\begin{lemma}
If $f(x)\geq 0 \ (\leq 0),\ x\in [\vp_1(T),\vp_1(T+U)]$ then

\be \label{9.6}
\begin{split}
& \int_T^{T+U}f[\vp_1(t)]\left|\zf\right|^2{\rm d}t= \\
& =\left\{ 1+\mcal{O}\left(\frac{\ln\ln T}{\ln T}\right)\right\}\ln T\int_{\vp_1(T)}^{\vp_1(T+U)}f(x){\rm d}x,\  U\in \left(\left. 0,\frac{T}{\ln T}\right.\right] ,
\end{split}
\ee

and

\be \label{9.7}
\begin{split}
& \int_{\mT}^{\widering{T+U}}f[\vp_1(t)]\left|\zf\right|^2{\rm d}t= \\
& =\left\{ 1+\mcal{O}\left(\frac{\ln\ln T}{\ln T}\right)\right\}\ln T\int_{T}^{T+U}f(x){\rm d}x,\ U\in \left(\left. 0,\frac{T}{\ln T}\right.\right] .
\end{split}
\ee
\end{lemma}

\begin{proof}
By using the mean-value theorem on the left-hand side of (\ref{9.3}) we directly obtain (\ref{9.6}), (see (\ref{9.2})). If we make use of the mean-value
theorem on the left-hand side of (\ref{9.5}) we obtain, by (\ref{9.2}),
\be \label{9.8}
\begin{split}
& \int_{\mT}^{\widering{T+U}}f[\vp_1(t)]\tilde{Z}^2(t){\rm d}t= \\
& \frac{1}{{\left\{ 1+\mcal{O}\left(\frac{\ln\ln t_1}{\ln t_1}\right)\right\}\ln t_1}}\int_{\mT}^{\widering{T+U}}
f[\vp_1(t)]\left|\zf\right|^2{\rm d}t ,
\end{split}
\ee

where $t_1\in (\mT,\widering{T+U})=(\vp_1^{-1}(T),\vp_1^{-1}(T+U))$ and

\be \label{9.9}
t_1=\vp_1^{-1}(T_1), T_1\in (T,T+U) .
\ee

Next, we obtain from (\ref{9.4}) by (\ref{9.9}) ($t_1\to\infty \ \Leftrightarrow \ T\to\infty$)

\be \label{9.10}
t_1-T_1=\mcal{O}\left(\frac{t_1}{\ln t_1}\right) \ \Rightarrow \ 1-\frac{T_1}{t_1}=\mcal{O}\left(\frac{1}{\ln t_1}\right)\to 0,\ T\to\infty ,
\ee

i.e.

\be \label{9.11}
t_1\sim T_1\sim T,\ T\to\infty ,
\ee

and (see (\ref{9.10}), (\ref{9.11}))

\be \label{9.12}
t_1-T=t_1-T_1+T_1-T=\mcal{O}\left(\frac{t_1}{\ln t_1}\right)+\mcal{O}(U)=\mcal{O}\left(\frac{T}{\ln T}\right) ,
\ee

where $U\leq\frac{T}{\ln T}$ by the condition of the Lemma 4. Now, (see (\ref{9.12}))

\be \label{9.13}
\ln t_1=\ln T+\mcal{O}\left(\frac{t_1-T}{T}\right)=\ln T+\mcal{O}\left(\frac{U}{T}\right)=\ln T+\mcal{O}\left(\frac{1}{\ln T}\right) .
\ee

Then the formula (\ref{9.7}) follows from (\ref{9.8}) by (\ref{9.1}), (\ref{9.13}).

\end{proof}

\section{Proof of Theorem 2}

\subsection{}

Putting

\bdis
f(t)=\left|\zf\right|^4
\edis

into (\ref{9.6}) we obtain

\be \label{10.1}
\int_T^{T+U_1}\left|\zfvp\right|^4\left|\zf\right|^2{\rm d}t\sim \ln T\int_{\vp_1(T)}^{\vp_1(T+U_1)}Z^4(t){\rm d}t ,
\ee

i.e. we have to consider the integral (see (\ref{8.2}))

\be \label{10.2}
\begin{split}
& \int_{\vp_1(T)}^{\vp_1(T+U_1)}Z^4(t){\rm d}t= \\
& \left\{\vp_1(t)\sum_{k=0}^4 C_k\ln^{4-k}\vp_1(t)\right\}_{t=T}^{t=T+U_1}+\mcal{O}(T^{7/8+\epsilon})= \\
& \sum_{k=0}^4C_kV_k+\mcal{O}(T^{7/8+\epsilon})
\end{split}
\ee

where

\be\label{10.3}
V_k=\vp_1(T+U_1)\ln^{4-k}\vp_1(T+U_1)-\vp_1(T)\ln^{4-k}\vp_1(T) ,
\ee

and, for example,

\be \label{10.4}
\begin{split}
& V_0=\vp_1(T+U_1)\ln^4\vp_1(T+U_1)-\vp_1(T)\ln^4\vp_1(T)= \\
& = \left. [\vp_1(T+U_1)-\vp_1(T)]\frac{{\rm d}}{{\rm d}\vp_1}[\vp_1(t)\ln^4(\vp_1(t))]\right|_{t=d_0} = \\
& = [\vp_1(T+U_1)-\vp_1(T)][\ln^4\vp_1(d_0)+4\ln^3\vp_1(d_0)] = \\
& =U_1\frac{\vp_1(T+U_1)-\vp_1(T)}{U_1}\ln^4\vp_1(d_0)\left\{ 1+\mcal{O}\left(\frac{1}{\ln\vp_1(d_0)}\right)\right\}, \\
& \vp_1(d_0)\in (\vp_1(T),\vp_1(T+U_1)) .
\end{split}
\ee

\subsection{}

By the Ingham formula (see \cite{3}, p. 294, \cite{9}, p. 120)

\bdis
\int_0^T Z^2(t){\rm d}t=T\ln T+(2c-1-\ln 2\pi)T+\mcal{O}(T^{1/2}\ln T)
\edis

we have in the case $U_1=T^{7/8+2\epsilon}$

\be \label{10.5}
\int_T^{T+U_1}Z^2(t){\rm d}t=U_1\ln T+(2c-\ln 2\pi)U_1+\mcal{O}(T^{7/8+\epsilon}) .
\ee

Comparing the formulae (\ref{2.5}) and (\ref{10.5}) we obtain

\be \label{10.6}
\frac{\vp_1(T+U_1)-\vp_1(T)}{U_1}=\tan[\alpha(T,U_1)]=1+\mcal{O}\left(\frac{\ln\ln T}{\ln T}\right) .
\ee

Since (see (\ref{10.6}))

\bdis
\vp_1(d_0)-\vp_1(T)\leq \vp_1(T+U_1)-\vp_1(T)=\mcal{O}(U_1) ,
\edis

we have (see (\ref{9.4}); $\vp_1(T)\sim T$)

\be \label{10.7}
\begin{split}
& \ln\vp_1(d_0)=\ln\vp_1(T)+\ln\left\{ 1+\frac{\vp_1(d_0)-\vp_1(T)}{\vp_1(T)}\right\}=\\
& =\ln\vp_1(T)+\mcal{O}\left(\frac{U_1}{T}\right)\sim \ln T .
\end{split}
\ee

Hence, we obtain from (\ref{10.4}) by (\ref{10.6}), (\ref{10.7})

\be \label{10.8}
V_0\sim U_1\ln^4 T ,
\ee

and similarly

\be\label{10.9}
V_l=\mcal{O}(U_1\ln^{4-l}T),\ l=1,2,3,4 .
\ee

Finally, the formula (\ref{8.5}) follows from (\ref{10.1}) by (\ref{10.2}), (\ref{10.8}) and (\ref{10.9}).

\section{Jacob's ladders and a new class of the nonlinear integral equations; concluding remarks}

\subsection{}

The proof of the Theorem 2 is simultaneously the proof of the following theorem.

\begin{theorem}
Every Jacob's ladder $\vp_1(t)=\frac 12\vp(t)$, where $\vp(t)$ is the exact solution of the nonlinear integral equation
\bdis
\int_0^{\mu[x(T)]}Z^2(t)e^{-\frac{2}{x(T)}t}{\rm d}t=\int_0^T Z^2(t){\rm d}t
\edis
is the asymptotic solution of the following nonlinear integral equation
\be \label{11.1}
\int_T^{T+U_1}\left|\zeta\left(\frac 12+ix(t)\right)\right|^4\left|\zf\right|^2{\rm d}t=\frac{1}{2\pi^2}U_1\ln^5T,\ U_1=T^{7/8+2\epsilon} ,
\ee
where $x(t)=x(t;T,\epsilon)$, for every fixed $T\geq T_0[\vp]$, i.e. the following asymptotic formula (see (\ref{8.5})
\bdis
\frac{2\pi^2}{U_1\ln^5T}\int_T^{T+U_1}\left|\zeta\left(\frac 12+i\vp_1(t)\right)\right|^4\left|\zf\right|^2{\rm d}t\sim 1, \ T\to\infty
\edis
holds true.
\end{theorem}

\subsection{}

Let us remind the Selberg's formula (\cite{8}, p. 128)

\be \label{11.2}
\int_T^{T+U_2}\{S(t)\}^{2k}{\rm d}\sim \frac{(2k)!}{k!(2\pi)^{2k}}U_2(\ln\ln T)^k ,
\ee

where $U_2=T^{1/2+\epsilon}$, and $k$ is the fixed positive number, and

\bdis
S(t)=\frac 1\pi\arg \zf
\edis

(where the $\arg$ is defined by the usual way). From (\ref{11.2}) by (\ref{9.7}) one obtains

\be \label{11.3}
\begin{split}
& \int_{\vp_1^{-1}(T)}^{\vp_1^{-1}(T+U_2)}\left\{\arg\zfvp\right\}^{2k}\left|\zf\right|^2{\rm d}t\sim  \\
& \sim \frac{(2k)!}{k!(2)^{2k}}U_2\ln T(\ln\ln T)^k,\ T\to\infty .
\end{split}
\ee

\begin{remark}
This formula cannot be obtained in the classical theory of A. Selberg and, all the less, in the theories of Balasubramanian, Heath-Brown and Ivic.
\end{remark}

Some nonlocal interaction of the functions

\bdis
\left\{\arg\zfvp\right\}^{2k},\ \left|\zf\right|^2
\edis

is expressed by the formula (\ref{11.3}).

\begin{remark}
Every Jacob's ladder $\vp_1(t)$ is the asymptotic solution (see (\ref{11.3})) of the nonlinear integral equation

\be \label{11.4}
\begin{split}
& \int_{x^{-1}(T)}^{x^{-1}(T+U_2)}\left\{\arg\zeta\left(\frac 12+ix(t)\right)\right\}^{2k}\left|\zf\right|^2{\rm d}t= \\
& =\frac{(2k)!}{k!2^{2k}}U_2\ln T(\ln\ln T)^k .
\end{split}
\ee
\end{remark}

\subsection{}

Since

\bdis
\int_T^{T+U}\pi(x){\rm d}x\sim \frac{UT}{\ln T},\ U\leq \frac{T}{\ln T}
\edis

holds true, we obtain (see (\ref{9.7}), $f(t)=\pi(t)$)

\be \label{11.5}
\int_{\vp_1^{-1}(T)}^{\vp_1^{-1}(T+U)}\pi[\vp_1(t)]\left|\zf\right|^2{\rm d}t\sim \frac{UT}{\ln T} ,
\ee

where $t-\vp_1(t)\sim (1-c)\pi(t)$ .

\begin{remark}
Every Jacob's ladder $\vp_1(t)$ is the asymptotic solution (see (\ref{11.5})) of the following nonlinear integral equation

\be \label{11.6}
\int_{x^{-1}(T)}^{x^{-1}(T+U)}\pi[x(t)]\left|\zf\right|^2{\rm d}t=\frac{UT}{\ln T},\ 0<U\leq \frac{T}{\ln T} .
\ee
\end{remark}

\subsection{}

Another source of the integrals containing the function $|\zf|^2$ is, for example, the system of the Chebyshev polynomials
$\mcal{T}_n(x),\ x\in [-1,1], n=0,1,2,\dots $
of the first kind. We obtain, from the well-known formula

\bdis
\int_{-1}^1\frac{[\mcal{T}_n(x)]^2}{\sqrt{1-x^2}}{\rm d}x=\left\{\begin{array}{rcl} \frac{\pi}{2} & , & n\geq 1 , \\ \pi & , & n=0 , \end{array} \right.
\edis

the following one

\bdis
\int_{T}^{T+2}\frac{[\mcal{T}_n(t-T-1)]^2}{\sqrt{1-(t-T-1)^2}}{\rm d}t=\frac{\pi}{2},\ n\geq 1 .
\edis

Next, putting

\bdis
f(t)=\frac{[\mcal{T}_n(t-T-1)]^2}{\sqrt{1-(t-T-1)^2}}
\edis

into (\ref{9.7}), $\mT=\vp_1^{-1}(T),\ \widering{T+2}=\vp_1^{-1}(T+2)$, we obtain

\be \label{11.7}
\int_{\vp_1^{-1}(T)}^{\vp_1^{-1}(T+2)}\frac{[\mcal{T}_n(\vp_1(t)-T-1)]^2}{\sqrt{1-(\vp_1(t)-T-1)^2}}\left|\zf\right|^2{\rm d}t\sim \frac{\pi}{2}\ln T,\ n\geq 1 .
\ee

\begin{remark}
Jacob's ladder $\vp_1(t)$ is the asymptotic solution of the nonlinear integral equation (see (\ref{11.2}))

\be \label{11.8}
\int_{x^{-1}(T)}^{x^{-1}(T+2)}\frac{[\mcal{T}_n(x(t)-T-1)]^2}{\sqrt{1-(\vp_1(t)-T-1)^2}}\left|\zf\right|^2{\rm d}t=\frac{\pi}{2}\ln T,\ n\geq 1 .
\ee
\end{remark}

\subsection{}

Let us remind the Liapunov equation

\be \label{11.9}
\gamma\int_{(a)}\frac{\rho'({\rm d}x')^3}{|\vec{r}-\vec{r}^\prime|}+\frac{1}{2}\omega_{ik}\omega_{jk}(x_i-a_i)(x_j-a_j)=V_a
\ee

(comp. \cite{fock}, pp. 334-337) for determining of the form of the integration domain $(a)$, (the density $\rho$ is prescribed), i.e. the equilibrium figures of the rotating
body.

Analogically to the case (\ref{11.9}), we will call the segment $[x^{-1}(T),x^{-1}(T+2)]$ entering the equation (\ref{11.8}), for example, the equilibrium segment and the
segment $[\vp_1^{-1}(T),\vp_1^{-1}(T+2)]$ will be called the asymptotical equilibrium segment.

\begin{remark}
By (\ref{11.7}) there is, for every fixed $T\geq T_0[\vp]$, a continuum of the asymptotic equilibrium segments $[\vp_1^{-1}(T),\vp_1^{-1}(T+2)]$. However, is there any
equilibrium segment $[x^{-1}(T),x^{-1}(T+2)]$ for some $T\geq T_0[\vp]$?
\end{remark}

\begin{remark}
There are the fixed-point methods and other methods of the functional analysis used to study the nonlinear equations. What can be obtained by using these
methods in the case of the nonlinear integral equations of the type (\ref{11.1}), (\ref{11.4}), (\ref{11.6}), (\ref{11.8})?
\end{remark}

\thanks{I would like to thank Ekatherina Karatsuba and Michal Demetrian for their moral support of my study of the Jacob's ladders.}

\end{document}